\documentclass[12pt,a4paper]{article}

\RequirePackage{etex}
\usepackage[utf8]{inputenc}
\usepackage{amsmath, amssymb, amsthm, amsfonts, mathrsfs,url}
\usepackage{bbm, commath}
\usepackage{multicol}
\usepackage{graphicx}
\usepackage{tikz-network}
\usetikzlibrary{positioning,patterns, calc}
\usetikzlibrary{graphs,graphs.standard,quotes}
\usepackage{authblk}
\usepackage{enumitem}

\usepackage{subcaption}
\usepackage{lipsum}
\usepackage{epstopdf}
\usepackage{algorithmic}
\usepackage[nameinlink]{cleveref}

\ifpdf
  \DeclareGraphicsExtensions{.eps,.pdf,.png,.jpg}
\else
  \DeclareGraphicsExtensions{.eps}
\fi

\crefname{hypothesis}{Hypothesis}{Hypotheses}
\newtheorem{thm}{Theorem}
\newtheorem{defn}{Definition}
\newtheorem{prop}{Proposition}
\newtheorem{lem}{Lemma}
\newtheorem{cor}{Corollary}

\newtheorem{rem}{Remark}

\def \Rl {{\mathbbm R}}
\def \Cl {{\mathbbm C}}
\def \e {{\bf e}}
\def \vl {{\bf v}}
\def \ul {{\bf u}}
\def \x {{\bf x}}
\def \y {{\bf y}}
\def \z {{\bf z}}
\def \w {{\bf w}}
\def \J {{\mathbf J}}
\def \l {{\mathbf 1}}
\def \o {{\bf 0}}

\newcommand{\ov}[1]{\left[#1\right]}
\newcommand{\ob}[1]{\left(#1\right)}
\newcommand{\cb}[1]{\left\lbrace #1\right\rbrace}
\newcommand{\tb}[1]{\left[#1\right]}
\newcommand{\mb}[1]{\left|#1\right|}
\newcommand{\up}[1]{\overset{#1}{\uplus}~}

\title{\LARGE{\it More graphs with pair state transfer}}

\author[1]{Hermie Monterde\thanks{Hermie.Monterde@uregina.ca}}
\affil[1]{\small{Department of Mathematics and Statistics, University of Regina, Regina, SK, Canada}}
\author[2]{Hiranmoy Pal\thanks{palh@nitrkl.ac.in}}
\affil[2]{\small{Department of Mathematics, NIT Rourkela, India-769008}}
\date{}
\begin{document}

\maketitle

% REQUIRED
\vspace*{-1 cm}
\begin{abstract}
This paper has two main goals. First, we characterize perfect state transfer between $s$-pair states in strongly regular graphs, as well as graphs in association schemes admitting perfect state transfer between vertices. The second goal is to provide a unified approach for constructing non-regular graphs admitting pair state transfer—relative to the adjacency, Laplacian, and signless Laplacian matrix—between the same pair of states at the same time. In particular, we show that for each $k\geq 5$, there are infinitely many connected graphs with maximum valency $k$ admitting this property. We also utilize graph products to generate new infinite families of graphs with pair state transfer.\\

% REQUIRED
{\bf keywords:}
quantum walk, perfect state transfer, graph spectra, adjacency matrix, Laplacian matrix, signless Laplacian matrix

% REQUIRED
\end{abstract}

\section{Introduction}
Let $G = (V, E, w)$ be a weighted graph with vertex set $V=V(G)$, edge set $E=E(G)$ and weight function $w : E \to \mathbb{R}^+$ that assigns a positive weight to each edge of $G$. We write $i\sim j$ if $A_{ij}\neq 0$. We say that $G$ is \emph{unweighted} if each edge weight in $G$ is one. The \emph{adjacency matrix} $A\in\mathbb{R}^{n\times n}$ of $G$ is defined by $A_{ij} = w(i,j)$ whenever $(i,j) \in E$, and $A_{ij} = 0$ otherwise. The \emph{degree matrix} $\Delta$ is the diagonal matrix such that $\Delta_{ii} = \sum_{j=1}^n A_{ij}$. The \emph{Laplacian matrix} and \emph{signless Laplacian matrix} of $G$ are given by $L = \Delta - A$ and $Q = \Delta + A$, respectively.

A \textit{(continuous) quantum walk} \cite{farhi} describes the time evolution of a quantum state on a quantum spin network, which is modeled by a graph $G$, where vertices represent qubits and edges correspond to interactions between them. When the coupling strengths between qubits are non-uniform, the underlying graph is weighted, with edge weights reflecting the interaction strengths. A quantum walk on a weighted graph $G$ defined by the \emph{transition matrix}
\begin{equation}\label{e1}
U_G(t)=\exp{\left(itM(G)\right)}=\sum\limits_{k\geq 0}\frac{(it)^k}{k!}M(G)^k,\quad t\in\Rl
\end{equation}
where $M(G)$ is taken to be the adjacency matrix $A$, Laplacian matrix $L$, or signless Laplacian matrix $Q$ of $G.$ If $G$ is regular, then $\Delta$ is a scalar matrix, and so the transition matrices governed by $A$, $L$ or $Q$ differ only by a global phase making the state transfer properties invariant under $A$, $L$, and $Q$. Unless stated otherwise, $M(G)$ denotes either $A$, $L$, or $Q$. If the context is clear, then we write $M=M(G)$. Throughout, $\J$ and $I$ resp.~denote the all-ones matrix and the identity matrix of appropriate sizes, $C_n$ and $CP(2m)$ resp.~denote the cycle and the cocktail party graph on $n$ and $2m$ vertices, and $Q_d$ denotes the hypercube of dimension $d$. We denote the complete bipartite graph with partite sets of sizes $m$ and $n$ by $K_{m,n}$. The complement of $G$ is denoted by $\overline{G}$. Vertices $a$ and $b$ in $G$ are \emph{twins} if they share the same neighbors. From \cite[Lem 2.9]{MonterdeELA}, $a$ and $b$ are twins in $G$ if and only if $(\e_a-\e_b)$ is an eigenvector for $A(G)$ associated with $\lambda=-1$ if $u\sim w$ and $\lambda=0$ otherwise.

\par The investigation of perfect state transfer was first initiated by the foundational contributions in \cite{bose, chr1}. In this context, a quantum state on a graph $G$ with $n$ vertices is represented by a unit vector in $\Cl^{n}$. The state corresponding to a vertex $a\in G$ is given by the characteristic vector $\e_a,$ while a superposition such as $\frac{1}{\sqrt{2}}\left(\e_a-\e_b\right)$ represents a pair state. More generally, quantum states are described by density matrices, which are positive semidefinite matrices with trace one. \emph{Perfect state transfer} (PST) occurs at time $\tau>0$ between two density matrices $D_1$ and $D_2$ whenever 
\begin{center}
$D_2=U_G(\tau)D_1U_G(-\tau).$
\end{center}
A \emph{real pure state} corresponding to a unit vector $\x\in\Rl^n$ is represented by the rank-one projection matrix $D_{\x}=\x\x^T$. Note that PST occurs between $D_\x$ and $D_y$ if and only if there exists $\tau>0$ and a \emph{phase factor} $\gamma\in\Cl$ such that
\begin{center}
$U_G(\tau)\x=\gamma\y.$
\end{center}
In this case, we simply say that PST occurs between $\x$ and $\y$. Since PST between vertex states is a rare phenomenon \cite{god2, god25}, a relaxation known as pretty good state transfer was introduced in \cite{god1}. This can also be extended to consider pretty good state transfer between real pure states. A graph $G$ admits \emph{pretty good state transfer} (PGST) between real pure states $\x$ and $\y$ if there exists a sequence $t_k\in\Rl$ and $\gamma\in\Cl$ such that
\begin{center}
$\displaystyle \lim\limits_{k\to\infty}U_G\ob{t_k}\x=\gamma\y.$
\end{center}
From \cite[Lem 22(3)]{god25}, PST between real pure states is \emph{monogamous}, i.e., if a real pure state $\x$ admits PST to both $\y$ and $\z$, then $\y = \z$. In contrast, PGST is not monogamous. As demonstrated in \cite[Ex 4.1]{pal5}, PGST can occur from a single vertex $u$ to multiple distinct vertices. The notion of a sedentary vertex was formally defined in\cite{mon}, which we extend to any $\x\in\Rl^n$. A real pure state $\x$ is \emph{$C$-sedentary} in $G$ if
\begin{center}
$\displaystyle\inf\limits_{t>0}~\mb{\x^TU_G(t)\x}\geq C$
\end{center}
for some constant $0 < C \leq 1$. As observed in \cite[Prop 2]{mon}, a $C$-sedentary state in a graph does not exhibit PGST.
\par Over the past two decades, a variety of graph families have been identified that admit PST and PGST. These include paths \cite{chr1, bom}, circulant graphs \cite{bas1, pal6}, Cayley graphs \cite{ber, che, pal2}, distance-regular graphs \cite{cou2}, and quotient graphs \cite{bac}. Additionally, vertex PST has been studied in graph products such as joins \cite{ange1,kirkland2026quantum}, coronas \cite{ack1}, non-complete extended p-sums (NEPS) \cite{li21, pal1}, and blow-up graphs \cite{pal9, mon2}. However, it has been shown that adjacency PST between vertices is rare \cite[Cor 10.2]{god1}. This lead to the introduction of PST between pair states, also known as \textit{pair state transfer}. Within the Laplacian framework, the concept of pair state transfer was initially introduced by Chen et al. \cite{che1}, where it was demonstrated that, among paths and cycles, only the paths on three or four vertices and the cycle on four vertices admit Laplacian pair state transfer, provided at least one pair forms an edge. Subsequently, it was shown in \cite{saro2} that any graph exhibiting vertex PST can serve as an isomorphic branch within a larger graph to facilitate pair state transfer. This construction lead to an infinite family of trees with maximum valency three admitting pair state transfer. Further developments concerning pair state transfer are presented in \cite{cao1, ojha2026laplacian,wang22}.
Recently, Kim et al. \cite{kim} generalized the notion of pair state transfer by developing the theory of PST between real $s$-pair states, which are quantum states of the form $\frac{1}{\sqrt{1+s^2}}\ob{\e_a+s\e_b}$, where $a\neq b$ and $s\in\Rl\backslash\{0\}$. A further generalization was presented in \cite{god25}, where PST between real pure states was investigated. For the sake of brevity, we omit the normalization factors $1/\sqrt{2}$ and $1/\sqrt{1+s^2}$ for pair states and $s$-pair states.

We organize this paper as follows. In section \ref{sec:as}, we characterize PST (resp., PGST) between $s$-pair states on graphs in association schemes admitting PST (resp., PGST) between vertices (\cref{mm1}). Since $CP(2m)$ for even $m\geq 2$ is the only family of strongly regular graphs admitting PST between vertices, this result motivates us in section \ref{sec:srg} to characterize PST between $s$-pair states in strongly regular graphs (\cref{imprimSRG1}). In section \ref{sec:main}, we investigate quantum state transfer in graphs with a structure that generalizes graphs with isomorphic branches and graphs with clusters 
(lemmas \ref{m1gen} and \ref{m1}). In section \ref{sec:pairst}, we use graphs with clusters to provide a unified approach for constructing non-regular graphs admitting pair state transfer—relative to $A,L,Q$—between the same pair of states at the same time (\cref{pairpst1}). We also show that for each $k\geq 5$, there are infinitely many connected graphs with maximum valency $k$ admitting this property (\cref{a}). We then apply our results to investigate pair state transfer in edge perturbed versions of complete graphs and cycles in section \ref{sec:compcyc}, and graph products in sections \ref{sec:seqj}, \ref{sec:comp} and \ref{sec:prods}. Finally, we present open questions in section \ref{sec:oq}. The remainder of this section is allotted to basic definitions and notation.

\par Let $\sigma(G)=\{\lambda_1,\ldots,\lambda_r\}$ be the set of all distinct eigenvalues of $M(G)$, and $E_{k}$ be the orthogonal projection onto the eigenspace associated with $\lambda_k\in \sigma(G)$. The spectral decomposition of the transition matrix is given by
\[U_G(t)=\sum\limits_{k=1}^r e^{it\lambda_k}E_{k}.\] The \textit{support} of a vector $\x$ in $G$ is the set $\sigma_{\x}(G)=\{\lambda_k\in\sigma(G):E_{k}\x\neq \o\}.$
A state $\x$ is \emph{fixed} if $\mb{\sigma_{\x}(H)}=1$, and this occurs precisely when $\x$ is an eigenvector for $M(H)$ \cite[Prop 2.3]{god25}. Fixed states do not admit PST \cite[Sec 4]{god25}. A state $\x$ is \textit{periodic} in $G$ whenever there is PST from $\x$ to itself. It is known that $\x$ is periodic if and only if $\sigma_{\x}(G)$ satisfies the \textit{ratio condition}, which is to say that $\frac{\lambda_j-\lambda_k}{\lambda_r-\lambda_s}$ for all $\lambda_j,\lambda_k,\lambda_r,\lambda_s\in \sigma_{\x}(G)$ with $r\neq s$. Two linearly independent states $\x$ and $\y$ in $\Rl^n$ are \textit{strongly cospectral} relative to $M(G)$ if $E_{k}\x=\pm E_{k}\y$
for all $\lambda_k\in\sigma_{\x}(G)$. In this case, $\sigma_{\x}(G)=\sigma_{\y}(G)$, and we may partition $\sigma_{\x}(G)$ into two nonempty sets given by
\begin{center}
$\sigma_{\x,\y}^+(G)=\{\lambda_k\in \sigma_{\x}(G):E_k\x=E_k\y\}\quad$ and $\quad \sigma_{\x,\y}^-(G)=\{\lambda_k\in \sigma_{\x}(G):E_k\x=-E_k\y\}$.
\end{center}
If $\x$ and $\y$ are strongly cospectral, then for each $\lambda\in\sigma_\x(G)$, either $\mathbf{v}^T\x=\mathbf{v}^T\y$ for any eigenvector $\mathbf{v}$ for $\lambda$ or $\mathbf{v}^T\x=-\mathbf{v}^T\y$ for any eigenvector $\mathbf{v}$ for $\lambda$. W

We close this section by stating the following characterizations of PST between real pure states \cite{god25}.

\begin{thm}
\label{PSTchar2}
Let $G$ be a graph on $n$ vertices and $\x,\y \in\Rl^n\backslash\{\mathbf{0}\}$ with $\y\neq \pm\x$. If $|\sigma_{\x}(G)|=2$, then $\x$ and $\y$ admit perfect state transfer if and only if they are strongly cospectral. In this case, the minimum PST time is $\frac{\pi}{|\lambda_1-\lambda_2|}$.
\end{thm}

\begin{thm}
\label{PSTchar}
Let $G$ be a graph on $n$ vertices and $\x,\y \in\Rl^n\backslash\{\mathbf{0}\}$ with $\y\neq \pm\x$. If $|\sigma_{\x}(G)|\geq 3$ and $\sigma_{\x}(G)$ is closed under algebraic conjugates, then $\x$ and $\y$ admit perfect state transfer if and only if all conditions below hold.
\begin{enumerate}
\item The vectors $\x$ and $\y$ are strongly cospectral.
\item Each eigenvalue $\lambda_j\in\sigma_{\x}(G)$ is of the form $\lambda_j=\frac{1}{2}(a+b_j\sqrt{\Delta})$, where $a$, $b_j$, and $\Delta$ are integers and either $\Delta=1$ or $\Delta>1$ is a square-free integer.
\item Let $\lambda_1\in\sigma_{\x,\y}^+(G)$. For all $\lambda_h\in\sigma_{\x,\y}^+(G)$ and $\lambda_{\ell},\lambda_k\in\sigma_{\x,\y}^-(G)$, $\nu_2\big(\frac{\lambda_1-\lambda_{h}}{\sqrt{\Delta}}\big)>\nu_2\big(\frac{\lambda_1-\lambda_{\ell}}{\sqrt{\Delta}}\big)=\nu_2\big(\frac{\lambda_1-\lambda_{k}}{\sqrt{\Delta}}\big).$
\end{enumerate}
In case there is PST between $\x$ and $\y$, the minimum PST time is $\frac{\pi}{g\sqrt{\Delta}}$, where $g=\operatorname{gcd}\big(\big\{\frac{\lambda_1-\lambda_j}{\sqrt{\Delta}}:\lambda_j\in\sigma_\x(G)\big\}\big)$.
\end{thm}

\section{Association schemes}
\label{sec:as}

A (symmetric) \emph{association scheme} with $d$ classes is a set $\mathcal{A}=\{A_0,...,A_d\}$ of symmetric 0-1 matrices such that $A_0=I$, $\sum_{j=0}^dA_j=\J$, and $A_iA_j$ belongs to the $\operatorname{span}\mathcal{A}$ for all $i,j$. The \textit{Bose-Mesner algebra} of an association scheme $\mathcal{A}$ is the complex algebra generated by $\mathcal{A}$. A \textit{graph belongs to an association scheme} if its adjacency matrix is contained in one. It is known that vertex-transitive graphs and distance-regular graphs belong to an association scheme. If $H$ is a connected graph in an association scheme that admits vertex PST at time $\tau$, then the $U_H(\tau)$ is a scalar multiple of a permutation matrix $P$ of order two without fixed points \cite[Sec 4]{cou2}. This fact extends to PGST: if vertex PGST occurs in $H$ with respect to a sequence $\{t_k\}$, then $\lim\limits_{k\to\infty}U_H(t_k)$ is a scalar multiple of a permutation matrix $P$ of order two without fixed points. We refer to the matrix $P$ as the permutation matrix associated with $H$.

\begin{thm}
\label{mm1}
Let $H$ be a connected graph in an association scheme that admits vertex PST (resp., PGST) with associated permutation matrix $P$. Then $H$ admits PST (resp., PGST) between $\x=\ob{\e_a+s\e_b}$ and $\y=\ob{\e_c+s\e_d}$ for all $s\in\Rl\backslash\{0\}$ if and only if either (i) $\e_c=P\e_a$ and $\e_d=P\e_b$, or (ii) $\e_b=P\e_a$ and $(c,d)=(b,a)$.
\end{thm}

\begin{proof}
Let $H$ be a connected graph in an association scheme that admits vertex PST (resp., PGST) between the pairs $\{a,c\}$ and $\{b,d\}$ at time $\tau$ with associated permutation matrix $P$.
First, suppose $P\e_a\neq \e_b$. Let $P\e_a=\e_c$ and $P\e_b=\e_d$. As $H$ admits vertex PST at time $\tau$, $U(\tau)\x=\alpha P\x=\alpha(P\e_a+sP\e_b)=\alpha\y$
for some $\alpha\in\Cl$. Next, suppose $P\e_a=\e_b$ and let $\x=\ob{\e_a+s\e_b}$. Then $P\e_b=\e_a$ and the same argument as the previous case yields $U(\tau)\x=\alpha\y$. The converse is immediate. The result for PGST follows similarly.  Finally, because PST is monogamous between real states \cite[Cor 5.3]{godsil2017real}, we get that $\x$ admits PST only with $\y$. For PGST, let $\z$ be a vector and suppose there is PGST between the pairs $\{\x,\y\}$ and $\{\x,\z\}$ relative to the sequences $\{\tau_k\}$ and $\{\tau_k'\}$, respectively. That is, $\lim\limits_{k\to\infty}U_H(\tau_k)\x=\gamma \y$ and $\lim\limits_{k\to\infty}U_H(\tau_k')\x=\gamma' \z$. Since $\lim\limits_{k\to\infty}U_H(\tau_k)=\alpha P$ and $\lim\limits_{k\to\infty}U_H(\tau_k')=\alpha'P$ for some $\alpha,\alpha'\in \Cl$, we get that $\y$ and $\z$ are scalar multiples of $P\x$. That is, $\x$ admits PGST only with $\y$.
\end{proof}

A \textit{Hadamard graph} of order $n\geq 2$ is a graph $G_\mathcal{H}$ obtained from an $n\times n$ Hadamard matrix $\mathcal{H}$ as follows: $G_\mathcal{H}$ has a pair of vertices $\{r^+,r^-\}$ for each row and a pair of vertices $\{c^+,c^-\}$ for each column. If the entry of $\mathcal{H}$ in row $r$ and column $c$ is 1, then $G_\mathcal{H}$ has edges $\{c^-,r^-\}$ and $\{c^+,r^+\}$. If the entry of $\mathcal{H}$ in row $r$ and column $c$ is $-1$, then $G_\mathcal{H}$ has edges $\{c^-,r^+\}$ and $\{c^+, r^-\}$. These graphs are distance-regular of diameter 4, antipodal with classes of size two, and bipartite with eigenvalues
\begin{center}
$\lambda_1=n,\quad$ $\lambda_2=\sqrt{n},\quad$ $\lambda_3=0,\quad$ $\lambda_4=-\lambda_1,\quad$ and $\quad\lambda_5=-\lambda_2$.
\end{center}
If $n$ is a perfect square, then it is known that Hadamard graphs of order $n$ admit vertex PST \cite[Cor 5.8]{cou2}. However, if $n$ is not a perfect square, then we show below that Hadamard graphs of order $n$ admit vertex PGST.

\begin{thm}
\label{hadpgst}
Let $G_{\mathcal{H}}$ be a Hadamard graph of order $n$. If $n$ is not a perfect square, then $G_{\mathcal{H}}$ admits vertex PGST. 
\end{thm}

\begin{proof}
From \cite[Section 5.2]{cou2}, we have that $E_j\e_u=E_j\e_v$ for $j\in\{1,3,5\}$ while $E_j\e_u=-E_j\e_v$ for $j\in\{2,4\}$. Since $\{n,\sqrt{n}\}$ is a linearly independent set, \cite[Lem 1]{kempton2017pretty} yields vertex PGST between any pair of antipodal vertices.
\end{proof}

If $H$ in \cref{mm1} is a distance-regular graph with vertex PST, then $H$ has $s$-pair PST for all $s\in\Rl\backslash\{0\}$. This recovers a result of Kim et al.\ \cite[Thm 6.6]{kim}. This applies if we take $H$ to be the hypercube $Q_d$ for any $d\geq 2$, a Hadamard graph with a perfect square order, or $CP(2m)$ with $m$ even. Moreover, if $H$ in \cref{mm1} is a Hadamard graph with a non perfect square order, then $H$ admits $s$-pair PGST for all $s\in\Rl\backslash\{0\}$ by \cref{hadpgst}.

\section{Strongly regular graphs}
\label{sec:srg}

A \textit{strongly regular graph} (SRG) with parameters $(n,k,a,c)$ is a $k$-regular graph on $n$ vertices such that (i) each pair of adjacent vertices has $a$ common neighbours and (ii) each pair of non-adjacent vertices has $c$ common neighbours. Here, we explicitly assume that $2\leq k\leq n-2$. For a strongly regular graph satisfying this condition, the eigenvalues of $A$ are $k$ and $\frac{1}{2}\left(a-c\pm \sqrt{(a-c)^2+4(k-c)}\right)$.
where we denote the positive eigenvalue by $\theta$ and the negative one by $\tau$. We denote the orthogonal projection matrices associated with $k$, $\theta$ and $\tau$ resp.~by $E_k=\frac{1}{n}J$, $E_\theta$ and $E_\tau$. A strongly regular graph $G$ is \textit{primitive} if $G$ and $\overline{G}$ are both connected. Otherwise, we say that $G$ is \textit{imprimitive}. It is known that a connected imprimitive strongly regular graph is isomorphic to the complete multipartite graph $\overline{\bigcup_m K_n}$ for some  $m,n\geq 2$ \cite[Section 1.1.3]{brouwer}.

The only family of connected strongly regular graphs with vertex PST is $CP(2m)$ for even $m\geq 2$. So by \cref{mm1}, this family also admits $s$-pair PST for all $s$. It is then natural to ask whether there are other strongly regular graphs that admit $s$-pair PST for some values of $s\neq 0$. In this section, we answer this question by characterizing PST between $s$-pair states in strongly regular graphs. 

\begin{prop}
\label{propsrg}
Let $G$ be a connected primitive strongly regular graph with $n\geq 5$ vertices. If $\x=\big(\e_u+s\e_v\big)$ and $\y=\big(\e_w+s\e_x\big)$ are strongly cospectral, then $E_\theta\x=\delta E_\theta\y$ and $E_\tau\x=-\delta E_\tau\y$ where $\delta\in\{\pm 1\}$.
\end{prop}

\begin{proof}
Since $E_k\x=E_k\y$, it cannot happen that $E_\theta\x=E_\theta\y$ and $E_\tau\x=E_\tau\y$. Otherwise, $\sigma_{\x,\y}^{-}(G)$ is empty, a contradiction. We now show that $E_\theta\x=-E_\theta\y$ and $E_\tau\x=-E_\tau\y$ also cannot happen. By way of contradiction, suppose $E_\theta\x=-E_\theta\y$ and $E_\tau\x=-E_\tau\y$. Then $(E_\theta+E_\tau)(\x+\y)=0$, i.e.
$(I-\frac{1}{n}J)(\e_u+\e_w+s(\e_v+\e_x))=\mathbf{0}$, which can write as $\frac{2(1+s)}{n}\mathbf{1}=\e_u+\e_w+s(\e_v+\e_x)$. Thus, $\frac{2(1+s)}{n}=1=s$, and so $n=4$, a contradiction.
\end{proof}

\begin{thm}
\label{srgchar}
If $G$ is a connected primitive strongly regular graph with $n\geq 5$ vertices, then $\x=(\e_u+s\e_v)$ and $\y=(\e_w+s\e_x)$ are not strongly cospectral in $G$ for all $u,v,w,x\in V(G)$ and for all $0\neq s\in\Rl$. Thus, PST between $s$-pair states does not occur in $G$.
\end{thm}

\begin{proof}
Fix $\mu\in\{\tau,\theta\}$ and let $\eta=\{\tau,\theta\}\backslash\{\mu\}$. As $G$ is primitive, $\theta>0$ and $\tau<-1$. Suppose $E_\mu\x=E_\mu\y$. Set $\z:=\x-\y$. As $E_\mu\z=E_k\z=0$ and $G$ only has three distinct eigenvalues, we get that $\z$ is an eigenvector for $\eta$, i.e. $A\z=\eta\z$. If $u=w$ or $v=x$, then $A\z=\eta\z$ implies that $(\e_v-\e_x)$ or $(\e_u-\e_w)$ are resp.~eigenvectors for $A$ associated with $\eta$. This implies that $\{v,x\}$ or $\{u,w\}$ is a twin set, in which case $\eta\in\{0,-1\}$, a contradiction. Thus, $u\neq w$ and $v\neq x$. Now, comparing corresponding entries in $A\z=\eta\z$ gives us
\vspace{-0.05in}
\begin{center}
\begin{multicols}{2}
$-A_{u,w}+s(A_{u,v}-A_{u,x})=\eta,$
$A_{u,w}+s(A_{w,v}-A_{w,x})=-\eta,$
\end{multicols}
\vspace{-0.3in}
\begin{multicols}{2}
$A_{u,v}-A_{v,w}=s(\eta+A_{v,x}),$
$A_{x,u}-A_{x,w}=-s(\eta+A_{v,x})$.
\end{multicols}
\end{center}
\vspace{-0.05in}
If $A_{u,v}=A_{u,x}$, then $\eta\in\{0,-1\}$ by the top left equation above, a contradiction. Thus $A_{u,v}\neq A_{u,x}$. Applying the same argument to other equations yields $A_{w,v}\neq A_{w,x}$, $A_{u,v}\neq A_{w,v}$ and $A_{x,u}\neq A_{x,w}$. As $A$ is a 0-1 matrix, we get  
\begin{equation}
\label{condii}
A_{u,v}=A_{x,w}\neq A_{u,x}=A_{w,v}.
\end{equation}
Moreover, $s=\pm (\eta+A_{u,w})$ and $s=\pm \frac{1}{\eta+A_{v,x}}$. Because $\theta>0$ and $\tau<-1$, the two preceding equations imply that $(\eta+A_{u,w})(\eta+A_{v,x})=1$. We have two cases.

\noindent \textbf{Case 1.} Suppose $G=C_5$ with vertices $0,1\ldots,5$ such that $j$ and $j+1$ are adjacent mod 5. Then $\theta=\frac{1}{2}(-1+\sqrt{5})$, $\tau=\frac{1}{2}(-1-\sqrt{5})$, and $
\mathbf{v}=(-1,-\theta,\theta,1,0)$ and $\mathbf{w}=(-1,-\tau,\tau,1,0)$ resp.~are eigenvectors associated with $\theta$ and $\tau$. In this case, $(\eta+A_{u,w})(\eta+A_{v,x})=1$ if and only if $A_{u,w}\neq A_{v,x}$. For the case when $A_{u,w}=1$, $A_{v,x}=0$ and $A_{u,v}=1$, the above observations above give us $A_{x,w}=1$, $A_{u,x}=0$, $s=\eta+1$. So, we may write $\x=(\e_0+s\e_4)$ and $\y=(\e_1+s\e_2)$. Since $\mathbf{v}$ is an eigenvector for $\theta$, $\mathbf{v}^T\x=-1$ and $\mathbf{v}^T\y=-\theta+(\eta+1)\theta$. This observation implies that if $\x$ and $\y$ are strongly cospectral, then $\eta=\tau$, and so $s=\tau
+1=-\theta$. Because $\mathbf{w}$ is an eigenvector for $\tau$, $\mathbf{w}^T\x=-1$ and $\mathbf{w}^T\y=-\tau-\theta\tau=-\tau+1\neq \pm 1$. Hence, $\x$ and $\y$ are not be strongly cospectral in $G$. The same argument applied to the cases (a) $A_{u,w}=1$, $A_{v,x}=0$ and $A_{u,v}=0$, (b) $A_{u,w}=0$, $A_{v,x}=1$ and $A_{u,v}=1$, and (c) $A_{u,w}=0$, $A_{v,x}=1$ and $A_{u,v}=0$ yields no strong cospectrality between $\x$ and $\y$ in $G$.

\noindent \textbf{Case 2.} Let $G\neq C_5$ so that $\eta\leq -2$. Then $(\eta+A_{u,w})(\eta+A_{v,x})=1$ if and only if $\eta=-2$ and
\begin{equation}
\label{condiii}
A_{u,w}=A_{v,x}=1.
\end{equation}
As $\theta>0$, we have $\eta=\tau$ and $\mu=\theta$, and thus $s=\pm 1$. We first deal with the case $s=1$. By \cref{propsrg}, $E_\theta\x= E_\theta\y$ and $E_\tau\x=- E_\tau\y$, which we may write as
$E_\theta(\e_u+\e_v)=E_\theta(\e_w+\e_x)$ and $ E_\tau\ul=\mathbf{0}$, where $\ul=(\e_u+\e_v+\e_w+\e_x)$. 
Now, observe that $\ul=(E_k+E_\theta+E_\tau)\ul=\frac{4}{n}\mathbf{1}+2E_\theta(\e_u+\e_v)$, and so we get
$E_\theta(\e_u+\e_v)=E_\theta(\e_w+\e_x)=\frac{1}{2}\ul-\frac{2}{n}\mathbf{1}$. As $(E_k+E_\theta+E_\tau)(\e_u+\e_v)=\e_u+\e_v$, we get $E_\tau(\e_u+\e_v)=\frac{1}{2}\ul$. Similarly, $E_\tau(\e_w+\e_x)=\frac{1}{2}\ul$, a contradiction to $E_\tau\x=- E_\tau\y$. Next, we examine the case $s=-1$. As $E_\tau\x=- E_\tau\y$, we have $E_\tau\vl=\mathbf{0}$, where $\vl=(\e_u-\e_v+\e_w-\e_x)$. Since we also have $E_k\vl=\mathbf{0}$, $\vl$ is an eigenvector associated with $\theta$. That is, $A\vl=\theta\vl,$
and thus  $-A_{u,v}+A_{u,w}-A_{u,x}=\theta$. Since $A_{u,w}=1$ from \cref{condiii}, we get $\theta=1-A_{u,v}-A_{u,x}$. By \cref{condii}, we have two subcases. If $A_{u,v}=1$, then $A_{u,x}=0$, and $\theta=0$, a contradiction. But if $A_{u,v}=0$, then $A_{u,x}=1$, and so $\theta=0$, a contradiction. Thus if $s=\pm 1$, then $\x$ and $\y$ are not strongly cospectral.

Combining cases 1 and 2, we get that strong cospectrality does not occur between $s$-pair states in $G$, and hence they do not admit PST by Theorems \ref{PSTchar2} and \ref{PSTchar}.
\end{proof}

\begin{thm}
\label{imprimSRG1}
Let $G$ be a connected strongly regular graph. Perfect state transfer occurs between $\x=\big(\e_u+s\e_v\big)$ and $\y=\big(\e_w+s\e_x\big)$ in $G$ if and only if one of the following holds.
\begin{enumerate}
\item $G=C_4$, between $(\e_u+\e_v)$ and $(\e_w+\e_x)$, where $u\not\sim v$ and $w\not\sim x$. 
\item $G=CP(2m)$, between (i) $(\e_u-\e_v)$ and $(\e_w-\e_x)$ for all $m\geq 2$, where $u\not\sim x$ and $v\not\sim w$ or (ii) $(\e_u+s\e_v)$ and $(\e_w+s\e_x)$ for all even $m\geq 2$, where $u\not\sim w$ and $v\not\sim x$. 
\item $G=\overline{\bigcup_m K_3}$ for all even $m$, between (i) $(\e_u+\frac{1}{2}\e_v)$ and $(\e_w+\frac{1}{2}\e_v)$ or (ii) $(2\e_u+\e_v)$ and $(2\e_w+\e_v)$, where $\{u,v,x\}$ belong to the same partite set. 
\item $G=\overline{\bigcup_m K_4}$ for all even $m$, between $(\e_u+\e_v)$ and $(\e_w+\e_x)$. 
\end{enumerate}
The minimum PST times in 1, 2, 3 and 4 are $\frac{\pi}{4}$, $\frac{\pi}{2}$, 
$\frac{\pi}{3}$, and $\frac{\pi}{4}$, respectively.
\end{thm}

\begin{proof}
By \cref{srgchar}, it suffices to focus on the case when $G$ is imprimitive, i.e. $G$ is isomorphic to $\overline{\bigcup_m K_n}$ for some  $m,n\geq 2$, which is $n(m-1)$-regular with $mn$ vertices and has eigenvalues $k=n(m-1)$, $\theta=0$ and $\tau=-n$, where $E_\theta$ is a direct sum of $m$ copies of the $n\times n$ matrix $I-\frac{1}{n}J_n$. 
In this case, $E_\theta\x=E_\theta\y$ if and only if 
\begin{equation}
\label{above}
\e_u-\e_w+s(\e_v-\e_x)=\frac{1}{n}\bigg(\sum_{j\not\sim u}\e_j-\sum_{j\not\sim w}\e_j+s\big(\sum_{j\not\sim v}\e_j-\sum_{j\not\sim x}\e_j\big)\bigg).
\end{equation}
Suppose \cref{above} holds. Note that the left-hand side (LHS) has four nonzero entries. If $u\not\sim w$ and $v\not\sim x$, then $\e_u+\e_w-\e_v-\e_x=\mathbf{0}$, a contradiction. If (i) $u\sim w$ and $v\not\sim x$ or (ii) $u\not\sim w$ and $v\sim x$, then the right-hand side (RHS) of \cref{above} has $2n$ nonzero entries. If $n=2$, then $\e_u-\e_w+s(\e_v-\e_x)$ equals $\frac{1}{2}\big(\sum_{j\not\sim u}\e_j-\sum_{j\not\sim w}\e_j\big)$ or $\frac{s}{2}\big(\sum_{j\not\sim v}\e_j-\sum_{j\not\sim x}\e_j\big)$. Both imply that $1=\frac{1}{2}$, a contradiction. Now, if $n\geq 3$, then the LHS of \cref{above} has less number of nonzero entries than the right, a contradiction. Finally, if $u\sim w$ and $v\sim x$, then we have two subcases. 
\begin{itemize}
\item Let $s\neq -1$. Then the RHS of \cref{above} has at least $2n$ nonzero entries. If it has more than $2n$ nonzero entries, then the LHS of \cref{above} has less number of nonzero entries than the right, a contradiction. If it is exactly $2n$, then $s=1$, $n=2$ $u\not\sim v$ and $w\not\sim x$. This implies that $E_\theta\x=0$. However, one checks that $E_\tau\x=-E_\tau\y$ if and only if $m=2$. Thus, $\x$ and $\y$ are strongly cospectral if and only if $m=n=2$. This yields $G=C_4$ which has PST between $(\e_u+\e_v)$ and $(\e_w+\e_x)$ \cite[Thm 6.5(iii)]{kim}. 

\item Let $s=-1$. Then $u\sim v$ and $w\sim x$ otherwise $\x$ or $\y$ is fixed. The same argument yields that the RHS of \cref{above} has exactly $2n$ entries. This occurs if and only if $n=2$, $u\not\sim x$ and $v\not\sim w$. In this case, $G=CP(2m)$ and we have $E_\tau\x=-E_\tau\y$. So $\sigma_{\x,\y}(G)=\{\theta,\tau\}$, and $\x=(\e_u-\e_v)$ and $\y=-(\e_x-\e_w)$ are strongly cospectral. Applying \cref{PSTchar2} yields PST between $\x$ and $\y$ with minimum time $\frac{\pi}{2}$.
\end{itemize}
The above cases prove 1 and 2(i). We now examine the case $E_\theta\x=-E_\theta\y$. Note that $E_\theta\x=-E_\theta\y$ if and only if 
\begin{equation}
\label{above1}
\e_u+\e_w+s(\e_v+\e_x)=\frac{1}{n}\bigg(\sum_{j\not\sim u}\e_j+\sum_{j\not\sim w}\e_j+s\big(\sum_{j\not\sim v}\e_j+\sum_{j\not\sim x}\e_j\big)\bigg).
\end{equation}
Suppose \cref{above1} holds. We start with case $s\neq -1$. We proceed with four subcases.
\begin{itemize}
\item Let $u,w,v,x$ belong to the same partite set. Then $\e_u+\e_w+s(\e_v+\e_x)=\frac{2(1+s)}{n}\sum_{j\not\sim u}\e_j$ implying that $s=1$, and $n=4$. In this case, $G$ is the complete multipartite graph $\overline{\bigcup_m K_4}$ and $E_\tau\x=-E_\tau\y$. So, $\x=(\e_u+\e_v)$ and $\y=(\e_w+\e_x)$ are strongly cospectral with $\sigma_{\x,\y}^+(G)=\{4(m-1),-4\}$ and $\sigma_{\x,\y}^-(G)=\{0\}$. By \cref{PSTchar}, we get PST between $\x$ and $\y$ if and only if $m$ is even, with minimum time $\frac{\pi}{4}$.
\item Suppose three vertices in $\{u,w,v,x\}$ belong to the same partite set. If $u$ is not one of these three vertices, then $\e_u+\e_w+s(\e_v+\e_x)=\frac{1}{n}\big(\sum_{j\not\sim u}\e_j+(2s+1)\sum_{j\not\sim v}\e_j\big)$ implying that $n=1$, a contradiction. The same argument for the case when $w$, $v$, or $x$ is not one of the three vertices yields the same contradiction.
\item Suppose two vertices in $\{u,v,w,x\}$ belong to a partite set that does not contain the other two. If the other two vertices are contained in different partite sets, then the RHS of \cref{above1} have $3n\geq 6$ nonzero entries, a contradiction. Thus, the other two vertices are contained in the same partite set. If $u\not\sim v$ and $w\not\sim x$, then we get the result in 1. The remaining case to consider is $u\not\sim w$ and $v\not\sim x$, which yields $\e_u+\e_w+s(\e_v+\e_x)=\frac{2}{n}\big(\sum_{j\not\sim u}\e_j+s\sum_{j\not\sim v}\e_j\big)$. This holds if and only if $n=2$. Thus, $G=CP(2m)$ and $E_\tau\x=E_\tau\y$. So $\x=(\e_u+s\e_v)$ and $\y=(\e_w+s\e_x)$ are strongly cospectral with $\sigma_{\x,\y}^+(G)=\{2(m-1),-2\}$ and $\sigma_{\x,\y}^-(G)=\{0\}$. By \cref{PSTchar}, we get PST between $\x$ and $\y$ if and only if $m$ is even, with minimum time $\frac{\pi}{2}$.
\item If $u,w,v,x$ are in distinct partite sets, then the RHS of \cref{above1} has $4n\geq 8$ nonzero entries, a contradiction.
\end{itemize}
The above subcases prove 2(ii) and 4. Next, let $s=-1$. Then $u\sim v$ and $w\sim x$, otherwise $\x$ or $\y$ is fixed. If $u\not\sim x$ and $w\not\sim v$, then $\e_u+\e_w-\e_v-\e_x=\mathbf{0}$, a contradiction. If $u\not\sim x$ and $w\sim v$, then $\e_u+\e_w-\e_v-\e_x=\frac{1}{n}\big(\sum_{j\not\sim w}\e_j-\sum_{j\not\sim v}\e_j\big)$ by \cref{above1}, and so $n=1$, a contradiction. Similarly for the case  $u\sim x$ and $w\not\sim v$. Finally, if $u\sim x$ and $w\sim v$, then we have two subcases. If $u\not\sim w$ and $v\not\sim x$, then \cref{above1}) gives us $\e_u+\e_w-\e_v-\e_x=\frac{2}{n}\big(\sum_{j\not\sim u}\e_j-\sum_{j\not\sim v}\e_j\big)$. This holds if and only if $n=2$, in which case we recover the result in 2(i) with the roles of $x$ and $w$ interchanged. However, if $u\sim w$ or $v\sim x$, then the RHS of \cref{above1} has $3n\geq 6$ nonzero entries, a contradiction.

To complete the proof, we examine the case when $u=w$ or $v=x$. Since PST between $(\e_u+s\e_v)$ and $(\e_u+s\e_x)$ is equivalent to PST between $s'\e_u+\e_v$ and $s'\e_u+\e_x$, where $s'=\frac{1}{s}$, we may assume without loss of generality that $u\neq w$ and $v=x$. Observe that $E_\theta\x=E_\theta\y$ if and only if $E_\theta(\e_u-\e_w)=\mathbf{0}$. Since $E_k(\e_u-\e_w)=\mathbf{0}$, it follows that $E_{\tau}(\e_u-\e_w)\neq \mathbf{0}$. That is, $(\e_u-\e_w)$ is an eigenvector associated with $\tau=-n$. Thus, $u$ and $w$ are twins, and so $\tau=0$ a contradiction. Thus, it cannot happen that $E_\theta\x=E_\theta\y$. Now, $E_\theta\x=-E_\theta\y$ if and only if
\begin{equation}
\label{above2}
\e_u+\e_w+2s\e_v=\frac{1}{n}\bigg(\sum_{j\not\sim u}\e_j+\sum_{j\not\sim w}\e_j+2s\sum_{j\not\sim v}\e_j\bigg).
\end{equation}
The LHS of \cref{above2} has three nonzero entries. If $u\sim w$ and $s\neq -\frac{1}{2}$, then the RHS of \cref{above2} has at least $2n\geq 4$ nonzero entries, a contradiction. If $s=-\frac{1}{2}$, then $\e_u+\e_w-\e_v=\big(\sum_{j\not\sim u}\e_j+\sum_{j\not\sim w}\e_j-\sum_{j\not\sim v}\e_j\big)$ implying that $n=1$, a contradiction. Now, suppose $u\not\sim w$. From \cref{above2}, $\e_u+\e_w+2s\e_v=\frac{2}{n}\big(\sum_{\not\sim u}\e_j+s\sum_{j\not\sim v}\e_j\big)$. If $u\sim v$, then the RHS has at least $2n\geq 4$ nonzero entries, a contradiction. If $u\not\sim v$, then $\e_u+\e_w+2s\e_v=\frac{2(1+s)}{n}\big(\sum_{j\not\sim u}\e_u\big)$ implying that $1=2s=\frac{2(1+s)}{n}$, i.e. $s=\frac{1}{2}$ and $n=3$. In this case, one checks that $E_\tau\x=E_\tau\y$. So $\x=(\e_u+\frac{1}{2}\e_v)$ and $\y=(\e_w+\frac{1}{2}\e_v)$ are strongly cospectral with $\sigma_{\x,\y}^+(G)=\{3(m-1),-3\}$ and $\sigma_{\x,\y}^-(G)=\{0\}$. By \cref{PSTchar}, PST occurs between $\x$ and $\y$ in $G$ if and only if $m$ is even, with minimum time $\frac{\pi}{3}$. So there is PST between $2\x=(2\e_u+\e_v)$ and $2\y=(2\e_w+\e_v)$ in $G$ at the same time. This proves 3.

The converse is straightforward.
\end{proof}

\section{Graph with special structure}
\label{sec:main}

Now we focus our attention to quantum walks on graphs with structure depicted in \cref{fgen}.
\begin{defn}\label{dgen}
Let $G$ be a connected graph with vertex set $B\cup C\cup S\cup T$ as in \cref{fgen} where the vertices in $B$ do not have neighbours in $S\cup T$, and the vertices in $C$ do not have neighbours in $T$.
\end{defn}
The matrix $M(G)$ is irreducible and may be written in the following block form
\begin{equation}\label{mgen}
    M(G)=\begin{bmatrix}
M[B,B] & M[B,C] & \o & \o \\
M[C,B] & M[C,C] & M[C,S] & \o \\
\o & M[S,C] & M[S,S] & M[S,T] \\
\o & \o & M[T,S] & M[T,T]
\end{bmatrix}.
\end{equation}
\begin{figure}
    \centering
\begin{tikzpicture}[yscale=0.8]
    \draw[thick] (0,0) ellipse (0.6cm and 1.2cm) node {$B$};
    \draw[thick] (2.5,0) ellipse (0.6cm and 1.2cm) node {$C$};
    \draw[thick] (5,0) ellipse (0.6cm and 1.2cm) node {$S$};
    \draw[thick] (7.5,0) ellipse (0.6cm and 1.2cm) node {$T$};
    
    \foreach \y in {-0.8,-0.4,0,0.4,0.8} {
        \draw[dashed] (0.6,\y) -- (1.9,\y);
        \draw[dashed] (3.1,\y) -- (4.4,\y);
        \draw[dashed] (5.6,\y) -- (6.9,\y);
    }
    
   %\draw[thick] (-0.3,-0.7) -- (2.8,-0.7);
   % \node at (1.25, -1.7) {$H$};
\end{tikzpicture}
\caption{A graph with special structure}
    \label{fgen}
\end{figure}

\begin{lem}\label{m1gen}
Suppose $G$ is a graph as in \cref{dgen}. Let $\widetilde{H}$ denote the subgraph of $G$ induced by $B\cup C$. Consider the original degrees in $G$ for the vertices of $\widetilde{H}$ to obtain
\begin{center}
    $M(\widetilde{H})=\begin{bmatrix}
M[B,B] & M[B,C] \\
M[C,B] & M[C,C]
\end{bmatrix},$
\end{center}
a principal submatrix of $M(G)$ as given in \eqref{mgen}. Suppose $U_{\widetilde{H}}(t)$ and $U_{G}(t)$ resp.~are the transition matrices of $\widetilde{H}$ and $G$. If $\x$ satisfies
\begin{equation}\label{gen}
\begin{bmatrix}\o & M[S,C]\end{bmatrix}M\left(\widetilde{H}\right)^{j}\x=\o, \quad \text{ for all } j\geq 0
\end{equation} then for all $t \in \mathbb{R}$,
\begin{center}
$U_{G}(t)\begin{bmatrix} \x\\
\o\end{bmatrix}=\begin{bmatrix} U_{\widetilde{H}}(t)\x\\ \o\end{bmatrix}.$
\end{center}
\end{lem}
\begin{proof}
Using the principle of mathematical induction, we observe that if $\x$ satisfies condition \eqref{gen}
then 
\begin{center}
$M(G)^{k}\begin{bmatrix} \x\\
\o\end{bmatrix}=\begin{bmatrix}M\left(\widetilde{H}\right)^{k}\x\\ \o\end{bmatrix}$
\end{center}
for all $k\ge0$. Hence the result follows from \cref{e1}.
\end{proof}

\subsection{Graphs with isomorphic branches}\label{subsec:isombr}

Let $X_1$ and $X_2$ be isomorphic graphs with isomorphism $f:V(X_1)\rightarrow V(X_2)$ satisfying $w(a,b)=w(f(a),f(b))$ for all $a,b\in V(X_1)$. Consider a connected graph $G$ where $X_1\cup X_2$, the disjoint union of $X_1$ and $X_2$, appears as an induced subgraph in such a way that $w(a,y)=w(f(a),y)$, for all $a\in V(X_1)$ and $y\in V(G)\backslash V(X_1\cup X_2)$. In this case, $X_1$ and $X_2$ are called \textit{isomorphic branches} in $G$. In \cite{saro2}, the authors investigated PST in graphs with isomorphic branches. We demonstrate that their main result \cite[Thm 20]{saro2} can be derived from \cref{m1gen}.

Let $\widetilde{H} = X_{1}\cup X_{2}$, where $X_{1}$ and $X_{2}$ are isomorphic branches of the graph $G$ under the mapping $f$. For $j=1,2$, we define $B_{j}=B\cap V(X_{j})$ and $C_{j}=C\cap V(X_{j})$, assuming that $f(B_1)=B_2$ and $f(C_1)=C_2$. Up to a reordering of the vertices, the adjacency matrix $A(G)$ takes the block form:
\[
A(G) = \left[
\begin{array}{cccc|cc}
A[B_1,B_1] & 0 & A[B_1,C_1] & 0 & 0 & 0 \\
0 & A[B_1,B_1] & 0 & A[B_1,C_1] & 0 & 0 \\
A[C_1,B_1] & 0 & A[C_1,C_1] & 0 & A[C_1,S] & 0 \\
0 & A[C_1,B_1] & 0 & A[C_1,C_1] & A[C_1,S] & 0 \\
\hline
0 & 0 & A[S,C_1] & A[S,C_1] & A[S,S] & A[S,T] \\
0 & 0 & 0 & 0 & A[T,S] & A[T,T]
\end{array}
\right].
\]
Here the adjacency matrix $A(\widetilde{H})$ of $\widetilde{H}$ appears as a principal submatrix of $A(G)$. Suppose $\y\in\Rl^{|B_1|}$ and $\z\in\Rl^{|C_1|}$. Letting $\x=[\y, -\y, \z, -\z]^T$, a simple induction shows that for all $k \ge 0$, the vector $A(\widetilde{H})^k \x$ takes the same form as that of $\x$. Consequently, condition \eqref{gen} of \cref{m1gen} is satisfied. Since $\widetilde{H}$ is the disjoint union of $X_1$ and $X_2$, if PGST occurs from a vertex $a$ in graph $X_1$ then pair PGST occurs from $(a,f(a))$ in $\widetilde{H}$. Hence, with a suitable choice of $\tilde{H}$, we obtain the conclusions established in \cite[Thm 2]{saro2} via \cref{m1gen}.

\subsection{Graphs with clusters}

Recall that vertices $a$ and $b$ in $G$ are \emph{twins} if they share the same neighbors. If, in addition, $a\not\sim b$, then they are called \emph{false twins}. Quantum state transfer between twin vertices was studied in \cite{kirk2,pal8}. In this section, we propose a generalized framework in which edges are added between a set of false twins $C$ in $G$. We call  $C$ a \textit{cluster} in $G$, first defined in \cite{merris}. 
The following definition extends the 
notion of a cluster in a graph.

\begin{defn}\label{df1}
Let $G$ be a graph and $C,S\subseteq V(G)$ with $|C|\geq 2$. The pair $(C,S)$ is a cluster in $G$ if $C$ has pairwise non-adjacent vertices, the neighborhood of each $a\in C$ in $G$ is equal to $S$, and for each $v\in S$, the edge weight $\omega(u,v)$ is constant for all $u \in C$.
\end{defn}

\begin{rem}
A connected graph with clusters is a special case of \cref{dgen} where $B = \emptyset$ and the row space of $M[S,C]$ in \eqref{mgen} is restricted to $\operatorname{span}\{\l^{T}\}$.
\end{rem}

\begin{defn}\label{df2}
Let $G$ be a connected graph with cluster $(C, S)$, and $H$ be a graph with vertex set $C$. We define $G(H)$ to be the graph with vertex set $V(G)$ and edge set $E(G) \cup E(H)$, where the edge weights from $G$ and $H$ are preserved.
\end{defn}

Let $G$ be a connected graph of order $n$ with cluster $(C, S)$ such that $|C|=c$ and $|S|=s$. Suppose the vertices in $C$ are labeled $1, 2, \ldots, c$, those in $S$ are labeled $c+1, c+2, \ldots, c+s$, and the remaining vertices are labeled $c+s+1,  \ldots, n$. Throughout, $H$ denotes a graph with vertex set $C$ such that $\l_c$ is an eigenvector of $M(H)$ (so that $H$ is regular whenever $M(H)\in\{A,Q\}$). As in \cref{df1}, for each $v \in S$, the edge weight $\omega(u,v)$ is constant for all $u \in C$. Let $\z \in \Rl^s$ be the vector whose components are $\z(v)=\omega(u,v),$ and define
\begin{center}
$\zeta=\begin{cases}-1 & \text{if $M=L$,}\\ 1 & \text{if $M\in\{A,Q\}$,} \end{cases}\quad \text{and}\quad \delta=\begin{cases}0 & \text{if $M=A$,}\\ 1 & \text{if $M\in\{L,Q\}$.} \end{cases}$
\end{center}
Let $\widetilde{M}(G-C)$ be the principal submatrix of $M(G)$ obtained by deleting all rows and columns indexed by the elements of $C.$ Following the natural ordering of the vertices of $G$, the matrices $M(G)$ and $M(G(H)$ are given by
\[
M(G) =
\left[
\begin{array}{cc}
\delta\, \l_s^T \z\, I_c & \hspace{0.6em} \tb{\zeta\, \l_c \z^T \quad \o} \\
[4pt]
\tb{\zeta\, \l_c \z^T \quad \o}^T & ~~\widetilde{M}(G - C)
\end{array}
\right]
\]
and
\[
M(G(H)) =
\left[
\begin{array}{cc}
\delta\, \l_s^T \z\, I_c + M(H) & \hspace{0.6em} \tb{\zeta\, \l_c \z^T \quad \o} \\
[4pt]
\tb{\zeta\, \l_c \z^T \quad \o}^T & ~~\widetilde{M}(G - C)
\end{array}
\right].
\]
Throughout, we let $\widetilde{\x}$ denote the vector $\ov{\x^T,\o}^T$.

\begin{lem}\label{ml1}
Let $\x\in\Rl^c$ satisfying $\l_c^T\x=0,$ and $H$ be a graph for which $\l_c$ is an eigenvector of $M(H)$. If $G(H)$ is the graph in \cref{df2}, then $\delta \l_s^T\z+\sigma_{\x}(H)=\sigma_{\widetilde{\x}}(G(H)).$
\end{lem}

\begin{proof}
If $\mathbf{1}_c$ is an eigenvector of $M(H)$, then the subspace of vectors $\mathbf{y}$ satisfying $\mathbf{1}_c^T\mathbf{y}=0$ can be decomposed into a direct sum of eigenspaces. So, any vector $\mathbf{y}$ orthogonal to $\mathbf{1}_c$ can be expressed as a linear combination of eigenvectors orthogonal to $\mathbf{1}_c$. Since $\mathbf{1}_c^T\mathbf{x}=0$, we get that $\mathbf{x}$ is an eigenvector for $H$ associated with the eigenvalue $\lambda$ if and only if $\widetilde{\mathbf{x}}$ is an eigenvector for $G(H)$ associated with $\delta \mathbf{1}_s^T\mathbf{z}+\lambda$. The desired result is then immediate.
\end{proof}

The following result establishes a relationship between the transition matrices of the graphs $H$ and $G(H).$ 

\begin{lem}\label{m1}
Suppose the conditions in \cref{ml1} are satisfied. If $U_H(t)$ and $U_{G(H)}(t)$ resp.~are the transition matrices of $H$ and $G(H)$, then for all $t \in \mathbb{R}$,
\begin{center}
$U_{G(H)}(t)\begin{bmatrix} \x\\
\o\end{bmatrix}=e^{i\delta t \l_s^T\z}\begin{bmatrix} U_H(t)\x\\ \o\end{bmatrix}.$
\end{center}
\end{lem}

\section{Pair state transfer}
\label{sec:pairst}

By \cref{m1}, the evolution of states of the form $\widetilde{\x}$ with $\l_c^T \x = 0$ in $G(H)$ depend solely on the subgraph $H$. The subspace $W \subseteq \Rl^c$ of vectors orthogonal to $\l_c$ is invariant under $M(H)$, as it admits a basis of eigenvectors of $M(H)$. Since $U_H(t)$ is a polynomial in $M(H)$, $W$ is also invariant under $U_H(t)$. So if $\x \in W$ and PST occurs between $\x$ and $\y$ in $H$, then $\y\in W$ and hence $\l_c^T \y = 0$. Likewise, if PST occurs in $G(H)$ from $\widetilde{\x}$ with $\l_c^T \x = 0$, then the transfer must be between states of the form $\widetilde{\x}$ and $\widetilde{\y}$ with $\l_c^T \y = 0$. The following is immediate from \cref{m1}.

\begin{thm}\label{mc1}
Let $G(H)$ be the graph in \cref{df2} and $\l_c$ be an eigenvector of $M(H)$ such that $\l_c^T \x =\l_c^T \y=0$.
\begin{enumerate}
    \item The state $\x$ is periodic (resp., sedentary) in $H$ if and only if $\widetilde{\x}$ is periodic (resp., sedentary) in $G(H)$.
    \item There is PST (resp., strong cospectrality, PGST) between $\x$ and $\y$ in $H$ if and only if there is PST (resp., strong cospectrality, PGST) between $\widetilde{\x}$ and $\widetilde{\y}$ in $G(H)$.
\end{enumerate}
Moreover, the time in which periodicity (resp., PST) occurs in $H$ and $G(H)$ in 1 (resp., 2) are the same. Similarly for the sequence of times for PGST in 2.
\end{thm}

Recall that the state transfer properties of a regular graph are invariant relative $A,L$ and $Q$. So if $H$ is regular with PST between $\x$ and $\y$ with $\l_c^T\x=0,$ then by \cref{mc1}, the graph $G(H),$ which may not be regular, admits PST between $\widetilde{\x}$ and $\widetilde{\y}$ irrespective of the above choices of $M(H).$ Consequently, if $H$ is regular, then \cref{mc1} holds for all $M\in\{A,L,Q\}$. Combining this fact with \cref{mc1} yields the following result.

\begin{figure}
\begin{center}
\begin{tikzpicture}[scale=.3,auto=left]
                       \tikzstyle{every node}=[circle, thick, fill=white, scale=0.5]
                       
		        \node[draw] (1) at (2,0) {u};
		        		        
		        \node[draw,minimum size=0.65cm, inner sep=0 pt] (2) at (-4, 2) {$a$};
		        \node[draw,minimum size=0.65cm, inner sep=0 pt] (3) at (-1, 2) {$c$};
		        \node[draw,minimum size=0.65cm, inner sep=0 pt] (4) at (-4, -2) {$b$};
		        \node[draw,minimum size=0.65cm, inner sep=0 pt] (5) at (-1, -2) {$d$};		        
                    
                \node[draw] (6) at (5,0) {};
                \node[draw] (7) at (8,0) {};
                \node[draw] (8) at (11,0) {};

				\draw [thick, black!70] (2)--(3)--(5)--(4)--(2);
				\draw [thick, black!70] (1)--(6)--(7);
				
				 \foreach \y in {2,3,4,5}{
                \draw[thick,black!70] (1)--(\y);}
                \foreach \x in {9,9.5,10} {\fill (\x,0) circle (2.5pt); }

                \node[draw] (w) at (27,0) {u};
                \node[draw] (c) at (22,0) {c};
                \node[draw,minimum size=0.65cm, inner sep=0 pt] (a) at (24.5, 2.5) {$a$};
		        \node[draw,minimum size=0.65cm, inner sep=0 pt] (b) at (24.5, -2.5) {$b$};	        
                    
                \node[draw] (f) at (30,0) {};
                \node[draw] (g) at (33,0) {};
                \node[draw] (h) at (36,0) {};

                \draw [thick, black!70] (w)--(a)--(c)--(b)--(w)--(f)--(g);
                \draw [thick, black!70] (w)--(c);
                \foreach \y in {a,b,c}{
                \draw[thick,black!70] (w)--(\y);}
                \foreach \x in {34,34.5,35} {\fill (\x,0) circle (2.5pt); }
		        \end{tikzpicture}
				
\end{center}	
\caption{\label{fig1} The graphs $G_1$ (left) and $G_2$ (right)}
\end{figure}

\begin{cor}
\label{pairpst1}
Let $G(H)$ as in \cref{df2} and $\l_c$ be an eigenvector of $M(H)$. Perfect state transfer occurs between $(\e_a-\e_b)$ and $(\e_c-\e_d)$ in $H$ if and only if it occurs between them in $G(H)$ at the same time. This holds for $A,L,Q$ whenever $H$ is a regular graph.
\end{cor}

Consider the graph $G_1$ in \cref{fig1}. The subgraph $H$ of $G_1$ induced by $\{a,b,c,d\}$ is $C_4$, which admits PST between $\ob{\e_a-\e_b}$ and $\ob{\e_c-\e_d}$ at $\frac{\pi}{2}$. As $\ob{V(H),\cb{u}}$ forms a cluster in $G_1-E(H),$ $G_1$ has PST between $\widetilde{\x}$ and $\widetilde{\y}$ at $\frac{\pi}{2}$ relative to $M\in\{A,L,Q\}$ by \cref{pairpst1}. This holds regardless of the length of the pendent path attached $u$. Thus:

\begin{thm}
\label{alqinf}
There are infinitely many connected non-regular unweighted graphs with maximum valency five that admit pair state transfer—relative to $A$, $L$, and $Q$—between the same pair states at $\frac{\pi}{2}$.
\end{thm}

Consider the graph $G_2$ in \cref{fig1}. The subgraph $H$ of $G_2$ induced by $\{a,b,c\}$ is the path on $3$ vertices, which admits Laplacian PST between $\ob{\e_a-\e_c}$ and $\ob{\e_b-\e_c}$ at $\frac{\pi}{2}$. The same argument above yields Laplacian PST between $\widetilde{\x}$ and $\widetilde{\y}$ in $G_2$ at $\frac{\pi}{2}$ by \cref{pairpst1}. This does not apply to $M\in\{A,Q\}$ as $\l_c$ is not an eigenvector for $M$.

\begin{thm}
There are infinitely many connected non-regular unweighted graphs with maximum valency four that admit Laplacian pair state transfer at $\frac{\pi}{2}$.
\end{thm}

\section{Complete graphs and cycles}
\label{sec:compcyc}

We show that the removal of a matching of size two from $K_n$ yields pair state transfer.

\begin{thm}\label{comp2}
Let $n\geq 4$. The removal of a matching of size at least two from $K_n$ results in pair state transfer relative to $A$, $L$ and $Q$ between the same pair of states at time $\frac{\pi}{2}$.
\end{thm}

\begin{proof}
If $\{a,b\}$ and $\{c,d\}$ are disjoint edges removed from $K_n$, then the resulting graph $G(C_4)$ can be viewed as a graph with cluster, where $C_4$ is the subgraph induced by $\{a,b,c,d\}$. In $\overline{G(C_4)}$, the edges $\{a,b\}$ and $\{c,d\}$ appear as a union of $K_2$. Since PST occurs in $K_2$ at $\frac{\pi}{2}$, we have PST between $(\e_a-\e_c)$ and $(\e_b-\e_d)$ in $K_2$ at the same time. Applying \cite[Thm 5.2]{che1}, we get PST between $(\e_a-\e_c)$ and $(\e_b-\e_d)$ in $G(C_4)$ relative to $L$. But since $C_4$ is regular, this result also applies to $A$ and $Q$ by \cref{pairpst1}.
\end{proof}

Let $K_n\backslash e$ denote the complete graph minus an edge $e=\{a,b\}$. By \cite[Cor 5.4]{che1}, $K_n\backslash e$ admits Laplacian PST between $(\e_a-\e_c)$ and $(\e_b-\e_c)$ for all $c\neq a,b$. The next result implies that two is the minimum size of a matching that needs to be removed from $K_n$ to obtain pair state transfer relative to $A,L,Q.$

\begin{prop}
For all $n\geq 4$, $K_n\backslash e$ does not admit pair state transfer relative to $A$.  
\end{prop}

\begin{proof} 
The eigenvalues of $A(K_n\backslash e)$ are $0$, $-1$, and $\lambda^{\pm}=\frac{1}{2}(n-3\pm\sqrt{(n+1)^2-8})$. The associated eigenvector for the eigenvalue $\lambda^+$, $\lambda^-$ and 0 resp.~are $\begin{bmatrix} -\lambda^-\l_2\\ 2\l_{n-2}\end{bmatrix}$, $\begin{bmatrix} -\lambda^+\l_2\\ 2\l_{n-2}\end{bmatrix}$, and $(\e_a-\e_b)$, where $\{a,b\}$ is the edge removed. Meanwhile, the set $\mathcal{B}=\{\e_c-\e_d:d\neq a,b,c\}$ forms a basis for the eigenspace associated with $-1$. Therefore, if $\x$ is a pair state that are is fixed, then it must be of the form $(\e_a-\e_c)$ and $(\e_b-\e_c)$, where $c\neq a,b$. In this case,  $\sigma_\x(G)=\{0,-1,\lambda^{\pm}=\frac{1}{2}(n-3\pm\sqrt{(n+1)^2-8})\}$, and so $\x$ is not periodic.
\end{proof}

In \cite[Thm 13]{pal4}, it is shown that the complement of $C_{2^k}$ admits PGST between antipodal vertices. The relative complement of $C_{2^k}$ in the complete graph $K_n$ with $n\geq 2^k$ vertices can be viewed as graph $G(H)$ with a cluster, where $H$ is the complement of $C_{2^k}$. The following is a direct consequence of \cref{mc1}(2).

\begin{thm}\label{comp3}
Let $n \geq 2^k$ and $C_{2^k}$ be a cycle embedded in $K_n$. The graph obtained by removing all edges of $C_{2^k}$ from $K_n$ has PGST between $(\e_a - \e_b)$ and $(\e_c - \e_d)$, where $\{a,c\}$ and $\{b,d\}$ are antipodal pairs of vertices in $C_{2^k}$.
\end{thm}

\section{Sequential join}
\label{sec:seqj}

The \emph{join} $G\vee H$ of graphs $G$ and $H$ is the graph obtained by taking the union of $G$ and $H$, and adding all possible edges between the vertices of $G$ and $H.$ 

\begin{cor}
\label{cor1}
Let $H$ be a connected graph in an association scheme that admits vertex PST (resp., PGST) with associated permutation matrix $P$. If  vertices $a,b\in V(H)$ do not admit PST, then $K_1\vee H$, $\overline{K_2}\vee H$, and $K_2\vee H$ admit PST (resp., PGST) between $\ob{\e_a-\e_b}$ and $\ob{\e_c-\e_d}$, where $\e_c=P\e_a$ and $\e_d=P\e_b$.
\end{cor}

\begin{proof}
We may view $K_1\vee H$ and $\overline{K_2}\vee H$ as $G(H)$, where $G\in\{K_{1,|V(H)|},K_{2,|V(H)|}\}$. If $G$ is $K_{2,|V(H)|}$ with an edge added in the partite set of two, then $K_2\vee H$ can be viewed as $G(H)$. The result then follows from \cref{cor1}.
\end{proof}
 
Suppose all edges in $G \vee H$ with one endpoint in $H$ incident to a fixed vertex in $G$ have identical weights. Then $\ob{V(H),V(G)}$ is a cluster in $(G\vee H)-E(H)$, in which case \cref{mc1} yields PST (resp., PGST) in $G\vee H$ between real states whenever $\l$ is an eigenvector for $H$ and $H$ admits PST (resp., PGST) between real states orthogonal to $\l$. This extends to the \emph{sequential join} $H_1 \vee H_2 \vee \cdots \vee H_k$, constructed by iteratively taking the join of each consecutive pair and forming the union of the resulting graphs. That is, $H_1 \vee H_2 \vee \cdots \vee H_k$ has vertex set $\bigcup_{j=1}^k V(H_j)$ and edge set $\bigcup_{j=1}^{k-1}
E(H_j\vee H_{j+1})$. For instance, the graphs in \cref{fig1} are $C_4 \vee K_1 \vee \cdots \vee K_1$ and $P_3 \vee K_1 \vee \cdots \vee K_1$.  

\begin{thm}\label{mc2}
$H_1 \vee H_2 \vee \cdots \vee H_k$ admits perfect (pretty good) real state transfer whenever some $H_j$ with eigenvector $\l$ admits perfect (pretty good) state transfer between two real states orthogonal to $\l$. In particular, $H_1 \vee H_2 \vee \cdots \vee H_k$ admits pair PST (resp., PGST) whenever some $H_j$ with eigenvector $\l$ admits pair PST (resp., PGST).
\end{thm}

Let $n\geq 4$ and $H_1$ be the graph obtained by removing a matching of size at least two from $K_n$. Then $H_1$ has pair state transfer relative to $A,L$ and $Q$ at $\frac{\pi}{2}$ by \cref{comp2}. Applying \cref{mc2}, $H_1\vee K_1\cdots\vee K_1$ inherits this property. Since this graph has maximum valency $n+1\geq 5$, we get an extension of \cref{alqinf} and \cite[Thm 5]{godsil2025quantum}.

\begin{thm}
\label{a}
For each $k\geq 5$, there are infinitely many connected non-regular unweighted graphs with maximum valency $k$ that admit pair state transfer —relative to $A$, $L$ or $Q$—between the same pair states at time $\frac{\pi}{2}$.
\end{thm}

By \cite[Theorem 3.8]{god25}, if we replace pair state transfer in \cref{a} with PST between pure states with nonnegative rational entries (e.g. vertex PST, PST between $(\e_u+\e_v)$ and $(\e_w+\e_x)$, etc.), then the result fails for all $k\geq 5$ relative to $M\in\{A,Q\}$. Nonetheless, there are infinite families of unweighted graphs that admit vertex PST—relative to $A$, $L$, and $Q$—between the same vertices, although these families do not have bounded valency and the PST times for $A,L,Q$ need not be the same. For example, $K_{2,4n+2}$ for all $n\geq 1$ admits PST between vertices in the partite set of size two at time $\tau=\frac{\pi}{2\sqrt{2n+1}}$ when $M=A$ and $\tau=\frac{\pi}{2}$ when $M\in\{L,Q\}$. This follows from \cite[Thm 12(1)]{kirk2}, \cite[Cor 5]{alv} and the fact that quantum walks on bipartite graphs are equivalent relative to $L$ and $Q$.

\section{Complement}
\label{sec:comp}

Let $\overline{M}=M(\overline{H})$. Note that $\overline{A} = \J - I - A$. If $H$ has $c$ vertices, then $\overline{L} =cI-\J-L$ and $\overline{Q} =(c-2)I+\J-Q$. Let $\l_c$ be an eigenvector of $M(H)$, and let $\x \in \mathbb{R}^c$ satisfy $\l_c^T \x = 0$. Since $\J$ commutes with $M(H)$ and $\J \x = 0$, the following relation holds between the transition matrices $U_H(t)$ and $U_{\overline{H}}(t)$ of $H$ and $\overline{H}$, respectively.

\begin{lem}\label{cg}
Let $\l_c$ be an eigenvector of $M(H)$, and let $\x \in \mathbb{R}^c$ satisfy $\l_c^T \x = 0$. Then $U_{\overline{H}}(t)\x = e^{i(\delta c - \zeta)t}~U_H(-t)\x$ for all $t\in\Rl$.
\end{lem}
 
The following is immediate from \cref{cg}.

\begin{thm}\label{cg1}
Let $H$ be a graph such that $\l_c$ be an eigenvector of $M(H)$. If $\x,\y \in \Rl^c$ with $\l_c^T \x =\l_c^T \y=0$, then:
\begin{enumerate}
    \item The state $\x$ is periodic (resp., sedentary) in $H$ if and only if $\x$ is periodic (resp., sedentary) in $\overline{H}$.
    \item There is PST (resp., strong cospectrality, PGST) between $\x$ and $\y$ in $H$ if and only if there is PST (resp., strong cospectrality, PGST) between $\x$ and $\y$ in $\overline{H}$.
\end{enumerate}
\end{thm}

We now combine \cref{cg1} with \cref{mc1} to derive further results. Consider the graph $G(H)$ in \cref{df2}. The complement $\overline{G(H)}$ corresponds to graph with cluster $\widetilde{G}(\overline{H})$, where $\widetilde{G}$ is obtained from $\overline{G}$ by removing all edges between the $c$ vertices that form the cluster in $G$. The following implications are now immediate.

\begin{thm}\label{cg2}
Suppose the premises in \cref{mc1} are satisfied. The following hold:
\begin{enumerate}
    \item The state $\x$ is periodic (resp., sedentary) in $H$ if and only if $\widetilde{\x}$ is periodic (resp., sedentary) in $\overline{G(H)}$.
    \item There is PST (resp., strong cospectrality, PGST) between $\x$ and $\y$ in $H$ if and only if there is PST (resp., strong cospectrality, PGST) between $\widetilde{\x}$ and $\widetilde{\y}$ in $\overline{G(H)}$.
\end{enumerate} 
\end{thm}
 
Observe that if $H$ is disconnected and a component of $H$ admits pair state transfer, then $H$ itself exhibits this property. The converse does not hold, since $K_2 \cup K_2$  admits pair state transfer, despite $K_2$ itself not exhibiting this property. In cases where $H$ in \cref{mc1} is disconnected, the following additional observations apply.

\begin{thm}\label{dmc1}
Let $H$ be a graph such that $\l_c$ is an eigenvector of $M(H)$ and $G(H)$ be the graph in \cref{df2}. If $H$ has two components, each admitting vertex PST (resp., PGST) at the same time with the same phase, then $G(H)$ and $\overline{G(H)}$ both admit pair PST (resp., pair PGST). In particular, if $H$ is regular, then this applies to $A,L$ and $Q$.
\end{thm}

\begin{figure}
\begin{center}
\begin{tikzpicture}[scale=.3,auto=left]
                       \tikzstyle{every node}=[circle, thick, fill=white, scale=0.4]
                       
		        \node[draw] (1) at (2,0) {};
		        		        
		        \node[draw] (2) at (-5, -1) {};
		        \node[draw] (3) at (-1, -2) {};
                
		        \node[draw] (4) at (-5, 1) {};
                \node[draw] (5) at (-1, 2) {};
		            
                         \node[draw] (6) at (5,0) {};
                         \node[draw] (7) at (8,0) {};
                         \node[draw] (8) at (11,0) {};

				\draw [thick, black!70] (2)--(3);
				\draw [thick, black!70] (1)--(6)--(7);
				
				\foreach \y in {2,3,4,5}{
                                   \draw[thick,black!70] (1)--(\y);
                                   }
                \foreach \x in {9,9.5,10} {\fill (\x,0) circle (2.5pt); }
                {\fill (-3.5,1.6) circle (2.5pt); }
                {\fill (-3,1.7) circle (2.5pt); }
                {\fill (-2.5,1.8) circle (2.5pt); }

		        	\end{tikzpicture}
				
\end{center}	
\caption{\label{fig7} A graph exhibiting Laplacian pair state transfer.}
\end{figure}

By \cref{dmc1}, adding a matching of size at least two within a partite set of $K_{m,n}$ yields Laplacian pair state transfer at $\frac{\pi}{2}$. Additionally, if $H$ is a perfect matching, then we get pair state transfer at $\frac{\pi}{2}$ relative to $A,L$ or $Q$.

\begin{thm}\label{dmc2}
Suppose the premises of \cref{dmc1} are satisfied. If $H$ has two components, one admitting vertex PST and the other periodic at the same time with the same phase, then $G(H)$ and $\overline{G(H)}$ admit pair state transfer.
\end{thm}

For any graph $H$ and for any $u\in V(H)$, we have $0\in \sigma_{\e_u}(H)$ relative to $L$. Thus, the phase factor for vertex PST relative to $L$ is always 1. Since $K_2$ admits vertex PST at $\frac{\pi}{2}$ and any isolated vertex can be regarded as Laplacian periodic at any time $t$, \cref{dmc2} implies that Laplacian pair state transfer at $\frac{\pi}{2}$ occurs in a graph formed by inserting an edge within a partite set of size at least three in $K_{m,n}$. This applies in general to graphs $G(H)$ with clusters, where $H$ contains an edge and an isolate vertex (see \cref{fig7}).

\section{More graph products} 
\label{sec:prods}

In this section, we construct graphs with pair PST using the Cartesian, corona, and lexicographic products.

The \emph{Cartesian product} $G_1 \square G_2$ of two graphs $G_1$ and $G_2$ is the graph with vertex set $V\ob{G_1} \times V\ob{G_2}$ satisfying
\begin{center}
$M\ob{G_1 \square G_2}=M\ob{G_1}\otimes I + I\otimes M\ob{G_2}$
\end{center}
for all $M\in\cb{A,L,Q}$. As $M\ob{G_1}\otimes I $ and $ I\otimes M\ob{G_2}$ commute, the transition matrix of $G_1 \square G_2$ can be evaluated as $U_{G_1 \square G_2}(t)=U_{G_1}(t)\otimes U_{G_2}(t).$ Consequently, the conclusions in \cite[Thm 47]{god25} hold when $M\in\cb{A,L,Q}$.  

\begin{thm}\label{cp1} 
     Let $H$ and $K$ be regular graphs and let $G(H)$ be the graph in \cref{df2}. Suppose $H$ admits PST between $\x$ and $\y$ at time $\tau$ and $\l_c^T \x=\l_c^T \y=0$. The following hold for all $M\in\cb{A,L,Q}.$
\begin{enumerate}
        \item If $\w$ is periodic in $K$ at time $\tau$, then $G(H)\square K$ exhibits PST between $\widetilde{\x}\otimes\w$ and $\widetilde{\y}\otimes\w$ at $\tau$.
        \item If $K$ admits PST between $\w_1$ and $\w_2$ at time $\tau$, then $G(H)\square K$ admits PST between $\widetilde{\x}\otimes\w_1$ and $\widetilde{\y}\otimes\w_2$ at $\tau$.
    \end{enumerate}
\end{thm}
From \cref{cp1}(2), if $H$ admits PST between pair states orthogonal to $\l_c$, and $K$ admits vertex PST (or has a periodic vertex), both at the same time, then $H\square K$ inherits pair state transfer relative to $A,L,Q$ between the same pair of states at the same time. For example, applying \cref{cp1}(2) to $(C_4 \vee K_1)\square K_2$ where $G(H)=C_4 \vee K_1$ and $H=C_4$ yields pair state transfer at $\frac{\pi}{2}$ when $M \in \{A, L, Q\}$. Note that $(C_4 \vee K_1)\square K_2$ is neither regular nor conforms to the structure of a graph with a cluster to which \cref{{pairpst1}} applies (see \cref{figcp}). Moreover, there are infinitely many such graphs of the form $(C_4 \vee K_1 \vee \cdots \vee K_1)\square K_2$. It is also worth noting that when $M=L$, the conclusions of \cref{cp1} hold even when $H$ is not regular. For example, $(K_2 \vee K_1 \vee K_1)\square K_2$ has Laplacian pair state transfer at $\frac{\pi}{2}$, owing to the fact that $H = K_2 \cup K_1$ itself admits this property.

\begin{figure}
\centering
\begin{tikzpicture}[scale=0.8,auto=left]
\tikzstyle{every node}=[circle, thick, fill=white, scale=0.4]
  \foreach \angle/\name in {45/1, 135/2, 225/3, 315/4}
    \node[draw] (a\name) at ({cos(\angle)}, {sin(\angle)}) {};

  \node[draw] (a0) at (0,0) {}; 

  \foreach \angle/\name in {45/1, 135/2, 225/3, 315/4}
    \node[draw] (b\name) at ({cos(\angle)+3}, {sin(\angle)+1}) {};

  \node[draw] (b0) at (3,1) {};

  \foreach \i/\j in {1/2, 2/3, 3/4, 4/1}
    \draw[thick] (a\i) -- (a\j);

  \foreach \i/\j in {1/2, 2/3, 3/4, 4/1}
    \draw[thick] (b\i) -- (b\j);
    
  \foreach \i in {1,2,3,4}
    \draw[thick] (a0) -- (a\i);

  \foreach \i in {1,2,3,4}
    \draw[thick] (b0) -- (b\i);

  \foreach \i in {0,1,2,3,4}
    \draw[thick] (a\i) -- (b\i);

\end{tikzpicture}
\caption{\label{figcp}The graph $P_2 \square (K_1 \vee C_4)$}
\end{figure} 

Now, let $G$ be a connected graph with $n$ vertices, and let $H$ be another graph with $c$ vertices. The \emph{vertex corona} $G \circ H$ of $G$ and $H$ is the graph obtained by taking one copy of $G$ and $n$ copies of $H$, and then joining the $i$-th vertex of $G$ to every vertex in the $i$-th copy of $H$ \cite{Barik2007}. The graph $G \circ H$ can be realized as a graph with multiple clusters. In the vertex corona $G \circ H$, labelling is carried out by first assigning labels to the vertices of $G$, followed by the labels of all vertices from the copies of $H$, ordered according to the vertices of $G$ to which they are attached. For any $a \in V(G)$ and $\x \in \mathbb{R}^c$, $\ob{\o,\, \e_a^T \otimes \x^T}$ is a real pure state in $\mathbb{R}^{n(1+c)}$. The following is immediate from \cref{mc1}.

\begin{thm}\label{corona1}
    Let $G$ be a graph with $a\in V(G)$, and let $\x \in \mathbb{R}^c$ with $\l_c^T \x = 0$. If $H$ is a graph such that $\l_c$ is an eigenvector of $M(H)$, then the following hold.
    \begin{enumerate}
    \item $\x$ is periodic (resp., sedentary) in $H$ if and only if $\ob{\o,\, \e_a^T \otimes \x^T}$ is periodic (resp., sedentary) in $G\circ H$.
    \item There is PST (resp., PGST) between $\x$ and $\y$ in $H$ if and only if there is PST (resp., PGST) between $\ob{\o,\, \e_a^T \otimes \x^T}$ and $\ob{\o,\, \e_a^T \otimes \y^T}$ in $G \circ H$. In particular, pair PST (resp., PGST) occurs in $H$ if and only if it occurs in $G \circ H$.
\end{enumerate}
\end{thm}

The above result generalizes  \cite[Thm 3.1]{wang22} on vertex corona. Analogous conclusions can be established for the edge and neighborhood corona \cite{Gopalapillai2011,Hou2010}, since these admit realizations as graphs with clusters, as illustrated in \cref{fig5}.

The \emph{blow-up} $\up{c}G$ of $c$ copies of a graph $G$ with $n$ vertices is the graph obtained by replacing each vertex of $G$ with an independent set of size $c$, such that vertices in different clusters are adjacent in $\up{c}G$ if and only if their corresponding vertices in $G$ are adjacent \cite{pal9}. Note that $\up{c}G$ consists of $n$ disjoint clusters, each with size $c$. If $c \geq 3$, then there is no vertex PST in $\up{c}G$ \cite[Cor 1(3)]{kirk2}. However, in such cases, additional edges can be introduced within those clusters to enable vertex PST relative to \cite{mon2}. Here we consider an edge-perturbed blow-up $\up{c}G\tb{H_1,H_2,\ldots,H_n}$ in which each cluster is replaced by a graph $H_j$ on $c$ vertices, for $j=1,2,\ldots,n$.

\begin{thm}\label{bl1}
Let $G$ be a graph on $n$ vertices. Then $\up{c}G\tb{H_1,H_2,\ldots,H_n}$ exhibits perfect (pretty good) real state transfer whenever some $H_j$ with eigenvector $\l$ admits perfect (pretty good) state transfer between two real states orthogonal to $\l$. In particular, $\up{c}G\tb{H_1,H_2,\ldots,H_n}$ admits pair PST (resp., PGST) whenever some $H_j$ with eigenvector $\l$ admits pair PST (resp., PGST).
\end{thm}

The \emph{lexicographic product} $G[H]$ of two graphs $G$ and $H$ is the graph with vertex set $V(G) \times V(H)$, where $(u, x)\sim(v, y)$ in $G[H]$ if either $u$ and $v$ are adjacent in $G$, or $u = v$ and $x\sim y$ in $H$. If each $H_j$ in $\up{c}G\tb{H_1,H_2,\ldots,H_n}$ is isomorphic to $H$, then the resulting graph is simply $G[H]$. We end this section with the following result.

\begin{figure}
\centering

% ===== Vertex Corona =====
\begin{minipage}{0.32\textwidth}
\centering
\begin{tikzpicture}[scale=.4,auto=left]
\tikzstyle{every node}=[circle, thick, fill=white, scale=0.4]

% Base P3 vertices
\node[draw] (u1) at (-2,0) {};
\node[draw] (u2) at (0,0) {};
\node[draw] (u3) at (2,0) {};

% Edges of C3
\draw[thick] (u1)--(u2)--(u3);

% C4 around u1
\node[draw] (a1) at (-3.5,2) {};
\node[draw] (b1) at (-0.5,2) {};
\node[draw] (c1) at (-1.5,1.5) {};
\node[draw] (d1) at (-2.5,1.5) {};

\draw[thick] (a1)--(b1)--(c1)--(d1)--(a1);
\foreach \v in {a1,b1,c1,d1} {
    \draw[thick] (u1)--(\v);
}

% C4 around u2
\node[draw] (a2) at (-1.5,-2) {};
\node[draw] (b2) at (1.5,-2) {};
\node[draw] (c2) at (0.5,-1.5) {};
\node[draw] (d2) at (-0.5,-1.5) {};

\draw[thick] (a2)--(b2)--(c2)--(d2)--(a2);
\foreach \v in {a2,b2,c2,d2} {
    \draw[thick] (u2)--(\v);
}
% C4 around u3
\node[draw] (a3) at (0.5,2) {};
\node[draw] (b3) at (3.5,2) {};
\node[draw] (c3) at (2.5,1.5) {};
\node[draw] (d3) at (1.5,1.5) {};

\draw[thick] (a3)--(b3)--(c3)--(d3)--(a3);
\foreach \v in {a3,b3,c3,d3} {
    \draw[thick] (u3)--(\v);
}
\end{tikzpicture}
\subcaption{Vertex corona}
\end{minipage}
%
% ===== Edge Corona =====
\begin{minipage}{0.32\textwidth}
\centering
\begin{tikzpicture}[scale=.4,auto=left]
\tikzstyle{every node}=[circle, thick, fill=white, scale=0.4]

% Base P3 vertices
\node[draw] (u1) at (-2,0) {};
\node[draw] (u2) at (0,0) {};
\node[draw] (u3) at (2,0) {};

% Edges of P3
\draw[thick] (u1)--(u2)--(u3);

% C4 for edge (u1, u2)
\node[draw] (e1a) at (-3,2) {};
\node[draw] (e1b) at (1,2) {};
\node[draw] (e1c) at (-0.5,1.5) {};
\node[draw] (e1d) at (-1.5,1.5) {};

\draw[thick] (e1a)--(e1b)--(e1c)--(e1d)--(e1a);
\foreach \v in {e1a,e1b,e1c,e1d} {
    \draw[thick] (u1)--(\v);
    \draw[thick] (u2)--(\v);
}

% C4 for edge (u2, u3)
\node[draw] (e2a) at (-1,-2) {};
\node[draw] (e2b) at (3,-2) {};
\node[draw] (e2c) at (1.5,-1.5) {};
\node[draw] (e2d) at (0.5,-1.5) {};

\draw[thick] (e2a)--(e2b)--(e2c)--(e2d)--(e2a);
\foreach \v in {e2a,e2b,e2c,e2d} {
    \draw[thick] (u2)--(\v);
    \draw[thick] (u3)--(\v);
}

\end{tikzpicture}
\subcaption{Edge corona}
\end{minipage}
%
% ===== Neighborhood Corona =====
\begin{minipage}{0.32\textwidth}
\centering
\begin{tikzpicture}[scale=.4,auto=left]
\tikzstyle{every node}=[circle, thick, fill=white, scale=0.4]

% Base P3 vertices
\node[draw] (u1) at (-1.5,0) {};
\node[draw] (u2) at (0,0) {};
\node[draw] (u3) at (1.5,0) {};

% Edges of P3
\draw[thick] (u1)--(u2)--(u3);

% C4 attached to u1: neighbor is u2
\node[draw] (n1a) at (-4,2) {};
\node[draw] (n1b) at (-0.5,2) {};
\node[draw] (n1c) at (-1,1.5) {};
\node[draw] (n1d) at (-2,1.5) {};

\draw[thick] (n1a)--(n1b)--(n1c)--(n1d)--(n1a);
\foreach \v in {n1a,n1b,n1c,n1d} {
    \draw[thick] (u2)--(\v);  % Only u2 is connected to this copy
}

% C4 attached to u2: neighbors are u1 and u3
\node[draw] (n2a) at (-3,-2) {};
\node[draw] (n2b) at (3,-2) {};
\node[draw] (n2c) at (0.5,-1.5) {};
\node[draw] (n2d) at (-0.5,-1.5) {};

\draw[thick] (n2a)--(n2b)--(n2c)--(n2d)--(n2a);
\foreach \v in {n2a,n2b,n2c,n2d} {
    \draw[thick] (u1)--(\v);  % u1 connects to this H copy
    \draw[thick] (u3)--(\v);  % u3 also connects to this H copy
}

% C4 attached to u3: neighbor is u2
\node[draw] (n3a) at (4,2) {};
\node[draw] (n3b) at (0.5,2) {};
\node[draw] (n3c) at (1,1.5) {};
\node[draw] (n3d) at (2,1.5) {};

\draw[thick] (n3a)--(n3b)--(n3c)--(n3d)--(n3a);
\foreach \v in {n3a,n3b,n3c,n3d} {
    \draw[thick] (u2)--(\v);  % Only u2 connects to this H copy
}

\end{tikzpicture}
\subcaption{Neighborhood corona}
\end{minipage}

\caption{\label{fig5} Variants of the corona of $P_3$ with $C_4$ exhibiting pair state transfer.}
\end{figure}

\begin{cor}\label{mc3}
$H_1 \vee H_2 \vee \cdots \vee H_k$ (resp., $G[H]$ and $G \circ H$) admits pair PST (resp., PGST) whenever at least one of the $H_j$'s (resp., $H$) is a connected graph that admits PST (resp., PGST) between two pairs of vertex states at the same time and with the same phase factor.
\end{cor}

\begin{proof}
The case when $H_1 \vee H_2 \vee \cdots \vee H_k$ follows from \cref{mm1} and \cref{mc2}, while that of $G[H]$ and $G \circ H$ resp.~follow from \cref{bl1} and \cref{corona1}(2).
\end{proof}

\section{Open questions}
\label{sec:oq}
We found two infinite families of strongly regular graphs that admit $s$-pair PST for some $s$ but do not admit vertex PST. Motivated by this result, we are interested in a characterization of distance-regular graphs and graphs in association schemes that admit $s$-pair PST but not vertex PST. We also ask, is there an infinite family of graphs that admit PST between states of the form $(\e_a-\e_b)$ and $(\e_c+\e_d)$? Since $\l$ is an eigenvector for $L$, the occurrence of PST between such states is not possible relative to $L$. Similarly, it cannot happen relative to $A$ when the graph is regular. 

\section*{Acknowledgments}
H.~Monterde is supported by the Pacific Institute for the Mathematical Sciences through the PIMS-Simons Postdoctoral Fellowship. H.~Pal gratefully acknowledges the support provided by the National Institute of Technology Rourkela, India. We are grateful to an anonymous reviewer for providing the work in Section 4 (specifically Lemma \ref{m1gen} and Section \ref{subsec:isombr}) that generalizes our result on graphs with clusters. 

%%%%%%% THE BIBLIOGRAPHY %%%%%%%
\bibliographystyle{abbrv}
\bibliography{References}

\end{document}